\documentclass[12pt,oneside,english]{amsart}
\usepackage[latin1]{inputenc}
\usepackage{a4wide}
\usepackage{amssymb}
\makeatletter

\newcommand{\mathcircumflex}[0]{\mbox{\^{}}}

 \theoremstyle{plain}    
 \newtheorem{thm}{Theorem}[section]
 \numberwithin{equation}{section} 
 \numberwithin{figure}{section} 
 \theoremstyle{plain}
 \theoremstyle{remark}    
 \newtheorem{notation}[thm]{Notation} 
 \theoremstyle{definition}
 \newtheorem{defn}[thm]{Definition}
 \theoremstyle{remark}
 \newtheorem{rem}[thm]{Remark}
 \theoremstyle{remark}    
 \newtheorem{claim}[thm]{Claim} 
 \newtheorem{observation}[thm]{Observation}
 \theoremstyle{definition}
 \newtheorem*{defn*}{Definition}
 \theoremstyle{plain}    
 \newtheorem{lem}[thm]{Lemma} 
 \newtheorem{main theorem}[thm]{Main Theorem} 
 \theoremstyle{definition}
 \newtheorem{example}[thm]{Example}
 \theoremstyle{remark}    
 \newtheorem{note}[thm]{Note} 
 \theoremstyle{remark}    
 \newtheorem{conclusion}[thm]{Conclusion} 
 \newtheorem{main conclusion}[thm]{Main Conclusion} 

\usepackage{babel}
\makeatother
\begin{document}

\title{The Automorphism Tower of a Centerless Group (Mostly) Without Choice}

\author{Itay Kaplan and Saharon Shelah}\thanks{The second author would like to thank the United States-Israel 
Binational Science Foundation for partial support of this research. Publication 882.
}

\begin{abstract}
For a centerless group $G$, we can define its automorphism tower. We define 
$G^{\alpha}$: $G^{0}=G$, $G^{\alpha+1}=Aut\left(G^{\alpha}\right)$
and for limit ordinals $G^{\delta}=\bigcup_{\alpha<\delta}G^{\alpha}$.
Let $\tau_{G}$ be the ordinal when the sequence stabilizes. Thomas' celebrated
theorem says $\tau_{G}<\left(2^{\left|G\right|}\right)^{+}$ and more.\\
If we consider Thomas' proof too set theoretical, we have here a shorter
proof with little set theory. However, set theoretically we get a
parallel theorem without the axiom of choice.\\
We attach to every element in $G^{\alpha}$, the $\alpha$-th member
of the automorphism tower of $G$, a unique quantifier free type over
$G$ (whish is a set of words from $G*\left\langle x\right\rangle $)
. This situation is generalized by defining {}``$\left(G,A\right)$
is a special pair''.
\end{abstract}
\maketitle

\section{Introduction}

\subsection*{background}

Given any centerless group $G$, we can embed $G$ into its automorphism
group $Aut(G)$. Since $Aut(G)$ is also without center, we can do
this again, and again. Thus we can define an increasing continuous
series $\left\langle G^{\alpha}\left|\alpha\in\mathbf{ord}\right.\right\rangle $
- The automorphism tower. The natural question that rises, is whether
this process stops, and when. We define $\tau_{G}=min\left\{ \alpha\left|G^{\alpha+1}=G^{\alpha}\right.\right\} $.\\
In 1939 (see \cite{key-3}) Weilant proved that for finite $G$, $\tau_{G}$
is finite. But there exist examples of centerless infinite groups
such that this process does not stop in any finite stage. For example
- the infinite dihedral group $D_{\infty}=\left\langle x,y\left|x^{2}=y^{2}=1\right.\right\rangle $
satisfies $Aut\left(D_{\infty}\right)\cong D_{\infty}$. So the question
remained open until 1984, when Simon Thomas' celebrated work (see
\cite{key-4}) proved that $\tau_{G}\leq\left(2^{\left|G\right|}\right)^{+}$.
He later (see \cite{key-2}) improved this to $\tau_{G}<\left(2^{\left|G\right|}\right)^{+}$.
\\
For a cardinal $\kappa$ we define $\tau_{\kappa}$ as the smallest
ordinal such that $\tau_{\kappa}>\tau_{G}$ for all centerless groups
$G$ of cardinality $\leq\kappa$. As an immediate conclusion from
Thomas' theorem we have $\tau_{\kappa}<\left(2^{\kappa}\right)^{+}$.
\\
We also define the normalizer tower of $H$ - a subgroup of a group
$G$ - in $G$: \\
$\left\langle nor_{G}^{\alpha}(H)\left|\alpha\in\mathbf{ord}\right.\right\rangle $
by $nor_{G}^{0}(H)=H$, $nor_{G}^{\alpha+1}(H)=nor(nor_{G}^{\alpha}(H))$
and $nor_{G}^{\delta}(H)=\bigcup\left\{ nor_{G}^{\alpha}(H)\left|\alpha<\delta\right.\right\} $
for $\delta$ limit. Let \\
$\tau_{G,H}=min\left\{ \alpha\left|nor_{G}^{\alpha+1}(H)=nor_{G}^{\alpha}(H)\right.\right\} $.
This construction turns out to be very useful, thanks to the following:\\
For a cardinal $\kappa$, let $\tau_{\kappa}^{nlg}$ be the smallest
ordinal such that $\tau_{\kappa}^{nlg}>\tau_{Aut\left(\mathfrak{A}\right),H}$,
for every structure of cardinality $\leq\kappa$ and $H\leq Aut\left(\mathfrak{A}\right)$
of cardinality $\leq\kappa$.\\
In \cite{key-1}, Just, Shelah and Thomas, found a connection between
these ordinals: $\tau_{\kappa}\geq\tau_{\kappa}^{nlg}$. \\
In this paper we deal with an upper bound of $\tau_{\kappa}$, but
there are conclusion regarding lower bounds as well, and the inequality
above is used to prove the existence of such lower bounds by finding
structures with long normalizer towers. In \cite{key-4}, Thomas proved
that $\tau_{\kappa}\geq\kappa^{+}$, and in \cite{key-1} the authors
found that one cannot prove in $ZFC$ a better explicit upper bound
for $\tau_{\kappa}$ then $\left(2^{\kappa}\right)^{+}$ (using set
theoretic forcing). In \cite{key-5}, Shelah proved that if $\kappa$
is strong limit singular of uncountable cofinality then $\tau_{\kappa}>2^{\kappa}$
(using results from $PCF$ theory).\\
It remains an open question whether or not there exists a countable
centerless group $G$ such that $\tau_{G}\geq\omega_{1}$.\\
In a subsequent paper we plan to prove that $\tau_{\kappa}^{nlg}\leq\tau_{\kappa}$ 
is true even without choice.

\subsection*{Results}

Our main theorem: (of course, Thomas' did not need to distinguish
$G$ and ${}^{\omega>}G$)

\begin{thm}
$\left(ZF\right)$ $\tau_{\left|G\right|}<\theta_{\mathcal{P}\left({}^{\omega>}G\right)}$
for a centerless group $G$. That is, there is an ordinal $\alpha$
and a function from ${}^{\omega>}G$ onto it such that $\tau_{G}<\alpha$.
Moreover, $\tau_{G'}<\alpha$ for every centerless group $G'$ such
that $\left|G'\right|\leq\left|G\right|$.
\end{thm}
This is an essentially theorem \ref{thm:ZFBoundOnTauK}.\\
We deal with finding $\tau_{G}$ without choice, and discover that
Thomas' theorem still holds. We prove that given a certain algebraic
property of $G$ and a subset $A$ ($\left(G,A\right)$ is special
- see definition \ref{def:Special}) we can reduce the bound. Along
the way we give a different proof of the theorem without choice in
conclusion \ref{con:Thomas} (Thomas used Fodor's lemma in his proof
as you can see in section \ref{sec:BackToChoice}, and it is known
that its negation is consistent with $ZF$) . Then we conclude that
if $\mathbf{V}'$ is a subclass of $\mathbf{V}$ which is a model
of $ZF$ such that $\mathcal{P}(\kappa)\in\mathbf{V}'$, then $\tau_{\kappa}<\left(\theta_{\mathcal{P}(\kappa)}\right)^{\mathbf{V}'}$
for every $\kappa$ and so $\tau_{\aleph_{0}}<\theta_{\mathbb{R}}^{L[\mathbb{R}]}$.
(see conclusion \ref{con:FianlConclusion})\\
Moreover, we give a descriptive set theoretic approach to finding
$\tau_{\aleph_{0}}$ in section \ref{sec:The-descriptive-set}. \\
Finally, we return to the axiom of choice, to see that we can improve
the bound for certain groups that hold a weaker algebraic property
($\left(G,A\right)$ is weakly special - see definition \ref{def:w-Special}).

\subsection*{A note about reading this paper}

How should you read this paper if you are not interested in the axiom
of choice but only on the new and simple proof of Thomas' Theorem?
\\
You can read only section \ref{sec:AutTowerWC}, and in there, you:
\\
Start with definition \ref{def:Special}. Continue to claim \ref{cla:MyLemma},
which is very simple. Then conclusion \ref{con:SpecialPairInd} is
a simple application of that lemma. Claim \ref{cla:FirstBound} Is
a very important step towards \ref{cla:secondBound}, and then finally
conclusion \ref{con:Thomas} wraps it up.

\begin{notation}
$ $
\begin{enumerate}
\item For a group $G$, its identity element, will be denoted as $e=e_{G}$.
\item if $A\subseteq G$ then$\left\langle A\right\rangle _{G}$ is the
subgroup generated by $A$ in $G$. Similarly, if $x\in G$ ,$\left\langle A,x\right\rangle _{G}$
is the subgroup generated by $A\cup\left\{ x\right\} $.
\item The language of a structure is its vocabulary.
\item $\mathbf{V}$ will denote the universe of sets; $\mathbf{V}'$ will
denote a transitive class which is a model of $ZF$.
\end{enumerate}
\end{notation}

\section{The Normalizer Tower Without Choice}

\begin{defn}
$ $
\begin{enumerate}
\item For a group $G$ and a subgroup $H\leq G$, we define $nor_{G}^{\alpha}\left(H\right)$
for every ordinal number $\alpha$ by:

\begin{itemize}
\item $nor_{G}^{0}\left(H\right)=H$.
\item $nor_{G}^{\alpha+1}\left(H\right)=nor_{G}\left(nor_{G}^{\alpha}\left(H\right)\right)$.
\item $nor_{G}^{\delta}\left(H\right)=\bigcup\left\{ nor_{G}^{\alpha}\left(H\right)\left|\alpha<\delta\right.\right\} $,
for $\delta$ limit.
\end{itemize}
\item We define $\tau_{G,H}^{nlg}=\tau_{G,H}=min\left\{ \alpha\left|nor_{G}^{\alpha+1}\left(H\right)=nor_{G}^{\alpha}\left(H\right)\right.\right\} $.
\item For a set $k$, we define $\tau_{\left|k\right|}^{nlg}$ as the smallest
ordinal $\alpha$, such that for every structure $\mathfrak{A}$ of
power $\left\Vert \mathfrak{A}\right\Vert \leq\left|k\right|$, $\tau_{Aut\left(\mathfrak{A}\right),H}<\alpha$
for every subgroup $H\leq Aut\left(\mathfrak{A}\right)=G$ of power
$\left|H\right|\leq\left|k\right|$. Note that $\tau_{\left|k\right|}^{nlg}=sup\left\{ \tau_{G,H}+1\left|\right.\right.$
for such $G$ and $\left.H\right\} $.
\item For a cardinal number $\kappa$, define $\tau_{\kappa}^{nlg}$ similarly.
\end{enumerate}
\end{defn}
\begin{rem}
Note that $\tau_{\left|k\right|}^{nlg}$ is well defined (in $ZF$)
since we can restrict ourselves to structures with languages of power
$\leq\sum_{n<\omega}\left|k\right|^{n}$ and universe contained in
$k$. See observation \ref{cla:LangKappa}.
\end{rem}
\begin{observation}$\,$
\label{cla:LangKappa}
\begin{enumerate}
\item ($ZF$) For any structure $\mathfrak{A}$ whose universe is $\left|\mathfrak{A}\right|=A$
there is a structure $\mathfrak{B}$ such that:

\begin{itemize}
\item $\mathfrak{A},\mathfrak{B}$ have the same universe (i.e. $A=\left|\mathfrak{B}\right|$).
\item $\mathfrak{A},\mathfrak{B}$ have the same automorphism group (i.e.
$Aut\left(\mathfrak{A}\right)=Aut\left(\mathfrak{B}\right)$).
\item the language of $\mathfrak{B}$ is of the form $L_{\mathfrak{B}}=\left\{ R_{\bar{a}}\left|\bar{a}\in{}^{\omega>}A\right.\right\} $
where each $R_{\bar{a}}$ is a $lg\left(\bar{a}\right)$ place relation.
\end{itemize}
\item ($ZFC$) If $\mathfrak{A}$ is infinite then the language of $\mathfrak{B}$
has cardinality at most $\left|A\right|$.
\end{enumerate}
\end{observation}
\begin{proof}
Define $\mathfrak{B}$ as follows: its universe is $\left|\mathfrak{A}\right|$.
Its language is $L=\left\{ R_{\bar{a}}\left|\bar{a}\in{}^{n}A,n<\omega\right.\right\} $
where $R_{\bar{a}}^{\mathfrak{B}}=o\left(\bar{a}\right)$, which is
defined by $o\left(\bar{a}\right)=\left\{ f\left(\bar{a}\right)\left|f\in Aut\left(\mathfrak{A}\right)\right.\right\} $
- the orbit of $\bar{a}$ under $Aut\left(\mathfrak{A}\right)$.
\end{proof}
\begin{defn}
For a set $A$, we define $\theta_{A}=\theta\left(A\right)$ to be
the first ordinal $\alpha>0$ such that there is no function from
$A$ onto $\alpha$.
\end{defn}
\begin{rem}
\label{rem:Theta}$\,$
\begin{enumerate}
\item $ZFC\vdash\theta_{A}=\left|A\right|^{+}$
\item $ZF\vdash$ $\theta_{A}$ is a cardinal number, and if $A$ is infinite
(i.e. there is an injection from $\omega$ into $A$) then $\theta_{A}>\aleph_{0}$.
\item Usually, we shall consider $\theta_{A}^{\mathbf{V}'}$ where $\mathbf{V}'$
is a transitive subclass of $\mathbf{V}$ which is a model of $ZF$.
\end{enumerate}
\end{rem}
\begin{claim}
\label{cla:ZFNormBound}$\left(ZF\right)$ If $G$ is a group, $H\leq G$
a subgroup then $\tau_{G,H}<\theta_{G}$.
\end{claim}
\begin{proof}
If $\tau_{G,H}=0$ it is clear. if not, define $F:G\to\tau_{G,H}$
by $F\left(g\right)=\alpha$ iff $g\in nor_{G}^{\alpha+1}\left(H\right)\backslash nor_{G}^{\alpha}\left(H\right)$,
and if there is no such $\alpha$, $F\left(g\right)=0$. By definition
of $\tau_{G,H}$, $F$ is onto. From the definition of $\theta$,
 $\tau_{G,H}<\theta_{G}$.
\end{proof}
We can do even more:

\begin{claim}
$\left(ZF\right)$ $\tau_{\left|k\right|}^{nlg}<\theta_{\mathcal{P}\left({}^{\omega>}k\right)}$.
\end{claim}
\begin{proof}
Let \begin{eqnarray*}
\mathcal{B}_{k} & = & \left\{ \left(\mathfrak{A},f,x\right)\right.\left|\textrm{$\mathfrak{A}$ is a structure, $L_{\mathfrak{A}}\subseteq{}^{\omega>}k$, $\left|\mathfrak{A}\right|\subseteq k$, $f:k\to{}^{k}k$, }\right.\\
 &  & \textrm{ $x\in G=Aut\left(\mathfrak{A}\right)$ and $\left.H\leq G,H=image\left(f\right)\right\} $}\end{eqnarray*}

Let $F:\mathcal{B}_{k}\to\tau_{\left|k\right|}^{nlg}$ be the following
map: $F\left(\mathfrak{A},f,x\right)=\alpha$ iff $x\in nor_{G}^{\alpha+1}\left(H\right)\backslash nor_{G}^{\alpha}\left(H\right)$,
and if there is no such $\alpha$, $F\left(G,H,x\right)=0$ (where
$G=Aut\left(\mathfrak{A}\right)$, and $H=image\left(f\right)$).
Since $F$ is onto $\tau_{\left|k\right|}^{nlg}$, its enough to show
that there is a one to one function from $\mathcal{B}_{k}$ to $\mathcal{P}\left({}^{\omega>}k\right)$.
But $x\in G$, hence $x\subseteq k\times k$ and $f\in{}^{k}\left({}^{k}k\right)\subseteq k\times k\times k$,
and $\mathfrak{A}$ is a series of subsets of ${}^{\omega>}k$, i.e.
a function in ${}^{{}^{\omega>}k}\mathcal{P}\left({}^{\omega>}k\right)\subseteq{}^{\omega>}k\times\mathcal{P}\left({}^{\omega>}k\right)$,
and we can encode such a function as a member of $\mathcal{P}\left({}^{\omega>}k\right)$.
(How? define an injective function $f_{1}:{}^{\omega>}k\times{}^{\omega>}k\to{}^{\omega>}k$,
using the definable injective function $cd:\omega\times\omega\to\omega$.
Then, define the encoding $f_{2}:{}^{\omega>}k\times\mathcal{P}\left({}^{\omega>}k\right)\to\mathcal{P}\left({}^{\omega>}k\right)$
using $f_{1}$). Hence it is clear.
\end{proof}
\begin{claim}
\label{cla:NormBound}Assume that $\mathbf{V}'$ is a transitive subclass
of $\mathbf{\mathbf{V}}$ which is a model of $ZF$, $G\in\mathbf{V}'$
a group, $H\in\mathbf{V}'$ a subgroup then $\tau_{G,H}^{\mathbf{V}}=\tau_{G,H}^{\mathbf{V}'}<\theta_{G}^{\mathbf{V}'}$.
\end{claim}
\begin{proof}
By claim \ref{cla:ZFNormBound}, it remains to show that $\tau_{G,H}^{\mathbf{V}}=\tau_{G,H}^{\mathbf{V}'}$.
By induction on $\alpha\in V'$, one can see that $\left(nor_{G}^{\alpha}\left(H\right)\right)^{\mathbf{V}}=\left(nor_{G}^{\alpha}\left(H\right)\right)^{\mathbf{V}'}$
(the formula that says that $x$ is in $nor_{G'}\left(H'\right)$
is bounded in the parameters $G'$ and $H'$). 
\end{proof}
It is also true that $\tau_{\left|k\right|}^{nlg}$ is preserved in
$\mathbf{V}'$, for every $k\in\mathbf{V}'$, such that $\mathcal{P}\left({}^{\omega>}k\right)\in\mathbf{V}'$:

\begin{claim}
\label{cla:NormTowerAbs}Assume that $\mathbf{V}'$ is a transitive
subclass of $\mathbf{V}$ which is a model of $ZF$.
\begin{enumerate}
\item \label{enu:NormAbsk}If $\mathcal{P}\left({}^{\omega>}k\right)\in\mathbf{V}'$
then $\left(\tau_{\left|k\right|}^{nlg}\right)^{\mathbf{V}'}=\left(\tau_{\left|k\right|}^{nlg}\right)^{\mathbf{V}}<\theta_{\mathcal{P}\left({}^{\omega>}k\right)}^{\mathbf{V}'}$.
\item \label{enu:NormAbsKappa}If $k=\kappa$ a cardinal number and $\mathcal{P}\left(\kappa\right)\in\mathbf{V}'$
then $\left(\tau_{\kappa}^{nlg}\right)^{\mathbf{V}'}=\left(\tau_{\kappa}^{nlg}\right)^{\mathbf{V}}<\theta_{\mathcal{P}\left(\kappa\right)}^{\mathbf{V}'}$
\end{enumerate}
\end{claim}
\begin{proof}
(\ref{enu:NormAbsKappa}) follows from (\ref{enu:NormAbsk}), as we
have an absolute definable bijection $cd:{}^{\omega>}\kappa\to\kappa$.\\
For a set $k\in\mathbf{V}'$, such that $\mathcal{P}\left({}^{\omega>}k\right)\in\mathbf{V}'$
let \begin{eqnarray*}
\mathcal{A}_{k} & = & \left\{ \left(G,H\right)\right.\left|\textrm{There is a structure $\mathfrak{A}$, with $\left|\mathfrak{A}\right|\subseteq k$, such that}\right.\\
 &  & \textrm{ $G=Aut\left(\mathfrak{A}\right)$ and $\left.H\leq G,\left|H\right|\leq\left|k\right|\right\} $}\end{eqnarray*}

It is enough to prove that $\left(\mathcal{A}_{k}\right)^{\mathbf{V}}=\left(\mathcal{A}_{k}\right)^{\mathbf{V}'}$,
because by definition\begin{eqnarray*}
\left(\tau_{\left|k\right|}^{nlg}\right)^{\mathbf{V}} & = & \bigcup\left\{ \tau_{G,H}+1\left|\left(G,H\right)\in\left(\mathcal{A}_{k}\right)^{\mathbf{V}}\right.\right\} \\
 & = & \bigcup\left\{ \tau_{G,H}+1\left|\left(G,H\right)\in\left(\mathcal{A}_{k}\right)^{\mathbf{V}'}\right.\right\} \\
 & = & \left(\tau_{\left|k\right|}^{nlg}\right)^{\mathbf{V}'}<\theta_{\mathcal{P}\left({}^{\omega>}k\right)}^{\mathbf{V}'}\end{eqnarray*}
 So let us prove the above equality: $\left(\mathcal{A}_{k}\right)^{\mathbf{V}'}\subseteq\left(\mathcal{A}_{k}\right)^{\mathbf{V}}$,
since if $\left(G,H\right)\in\left(\mathcal{A}_{k}\right)^{\mathbf{V}'}$
and $\mathfrak{A}\in\mathbf{V}'$ a structure such that $G=Aut\left(\mathfrak{A}\right)$
then $\mathfrak{A}\in\mathbf{V}$ and $\left(Aut\left(\mathfrak{A}\right)\right)^{\mathbf{V}}=\left(Aut\left(\mathfrak{A}\right)\right)^{\mathbf{V}'}$,
because $\left(Aut\left(\mathfrak{A}\right)\right)^{\mathbf{V}}\subseteq{}^{k}k\subseteq\mathcal{P}\left(k\times k\right)\in\mathbf{V}'$.
So $\left(G,H\right)\in\left(\mathcal{A}_{k}\right)^{\mathbf{V}'}$,
as witnessed by the same structure.\\
On the other hand, suppose $\left(G,H\right)\in\left(\mathcal{A}_{k}\right)^{\mathbf{V}}.$
So let $\mathfrak{A}$ be a structure on $k$ such that $G=Aut\left(\mathfrak{A}\right)$.
By observation \ref{cla:LangKappa}, we may assume that $L_{\mathfrak{A}}=\left\{ R_{\bar{a}}\left|\bar{a}\in{}^{\omega>}k\right.\right\} $,
and each $R_{\bar{a}}$ is a $lg\left(\bar{a}\right)$ place relation
(This is not necessary, it just makes it more convenient). Define
$X_{\mathfrak{A}}=\left\{ \bar{a}\mathcircumflex\bar{b}\left|lg\left(\bar{a}\right)=lg\left(\bar{b}\right)\land\bar{b}\in R_{\bar{a}}^{\mathfrak{A}}\right.\right\} $.
Observe that:
\begin{itemize}
\item $X_{\mathfrak{A}}\in\mathbf{V}'$, as $X_{\mathfrak{A}}\subseteq{}^{\omega>}k$.
\item $\mathfrak{A}$ can be defined using $X_{\mathfrak{A}}$: its universe
is $k$, and for each $\bar{a}\in{}^{\omega>}k$, $R_{\bar{a}}=\left\{ \bar{b}\left|lg\left(\bar{b}\right)=lg\left(\bar{a}\right)\land\bar{a}\mathcircumflex\bar{b}\in X_{\mathfrak{A}}\right.\right\} $.
\end{itemize}
So in conclusion, $\mathfrak{A}\in\mathbf{V}'$, and so $G\in\mathbf{V}'$
as before. In addition $H\in\mathbf{V}'$, because $H$ is the image
of a function in ${}^{k}\left({}^{k}k\right)$, and ${}^{k}\left({}^{k}k\right)\subseteq\mathcal{P}\left(k\times k\times k\right)\in\mathbf{V}'$.
By definition $\left(G,H\right)\in\left(\mathcal{A}_{k}\right)^{\mathbf{V}'}$
and we are done.
\end{proof}

\section{\label{sec:AutTowerWC}The Automorphism Tower Without Choice}

\begin{defn}
For a centerless group $G$, we define the series $\left\langle G^{\alpha}\left|\alpha\in\mathbf{ord}\right.\right\rangle $:
\end{defn}
\begin{itemize}
\item $G^{0}=G$.
\item $G^{\alpha+1}=Aut\left(G^{\alpha}\right)$
\item $G^{\delta}=\cup\left\{ G^{\alpha}\left|\alpha<\delta\right.\right\} $
for $\delta$ limit.
\end{itemize}
\begin{rem}
Since $G$ is centerless, this makes sense - $G\cong Inn\left(G\right)\leq Aut\left(G\right)$,
and $Aut\left(G\right)$ is again without center. So we identify $G$
with $Inn\left(G\right)$, and so $G^{\alpha}\leq G^{\alpha+1}$.
This series is therefore monotone and continuous.
\end{rem}
\begin{defn}
$\,$
\end{defn}
\begin{enumerate}
\item Define an ordinal $\tau_{G}$ by $\tau_{G}=min\left\{ \alpha\left|G^{\alpha+1}=G^{\alpha}\right.\right\} $.
We shall show below that $\tau_{G}$ is well defined.
\item For a set $k$, we define $\tau_{\left|k\right|}$ to be the smallest
ordinal $\alpha$ such that $\alpha>\tau_{G}$ for all groups $G$
with power $\leq\left|k\right|$.
\item For a cardinal number $\kappa$, define $\tau_{\kappa}$ similarly.
\end{enumerate}
\begin{defn}
For a group $G$ (not necessarily centerless) and a subset $A$, we
define an equivalence relation $E_{G,A}$ by $xE_{G,A}y$ iff $tp_{qf}\left(x,A,G\right)=tp_{qf}\left(y,A,G\right)$
where $tp_{qf}\left(x,A,G\right)=$ \begin{eqnarray*}
\left\{ \sigma\left(z,\bar{a}\right)\left|\right.\right. & \bar{a}\in{}^{n}A,n<\omega,\sigma\textrm{ a term in the language of groups (i.e. a word)}\\
 & \textrm{with parameters from }A,\\
 & \left.z\textrm{ is it's only free variable and }G\models\sigma\left(x,\bar{a}\right)=e\right\} \end{eqnarray*}

\end{defn}
\begin{rem}
\label{rem:E_GA}$\,$
\end{rem}
\begin{enumerate}
\item \label{enu:EqivDefOfE_GA}Note that $xE_{G,A}y$ iff there is an isomorphism
between $\left\langle A,x\right\rangle _{G}$ and $\left\langle A,y\right\rangle _{G}$
taking $x$ to $y$ and fixing $A$.
\item The relation $E_{G,A}$ is definable and absolute (since $tp_{qf}\left(x,A,G\right)$
is absolute - the formula defining it is bounded).
\end{enumerate}
\begin{defn}
\label{def:Special}We say $\left(G,A\right)$ is a special pair if
$A\subseteq G$, $G$ is a group and $E_{G,A}=\left\{ \left(x,x\right)\left|x\in G\right.\right\} $
(i.e. the equality).
\end{defn}
\begin{example}
$\,$
\end{example}
\begin{enumerate}
\item \label{ACreatesG}If $G=\left\langle A\right\rangle _{G}$ then $\left(G,A\right)$
is special.
\item \label{G_AutGIsSpecial}If $G$ is centerless then $\left(Aut\left(G\right),G\right)$
is special (see claim \ref{cla:MyLemma}), so in general, the converse
of (\ref{ACreatesG}) is not true.
\item There is a group $G$ with center such that $\left(Aut\left(G\right),Inn\left(G\right)\right)$
is special, e.g. $\mathbb{Z}/2\mathbb{Z}$, but
\item If $G$ is not centerless then (\ref{G_AutGIsSpecial}) is not necessarily
true, even if the $\left|Z\left(G\right)\right|=2$:

It is enough to find a group which satisfies these properties:

\begin{enumerate}
\item $Z\left(G\right)=\left\{ a,e\right\} $ where $a\neq e$.
\item $H_{i}\leq G$ for $i=1,2$ are two different subgroups of index 2.
\item $Z\left(G\right)=Z\left(H_{i}\right)$ for $i=1,2$
\end{enumerate}
Let $\pi$ be the homomorphism $\pi:G\to Aut\left(G\right)$ taking
$g$ to $i_{g}$ ($i_{g}\left(x\right)=gxg^{-1}$). Then $inn\left(G\right)=image\left(\pi\right)$.
We wish to find $x_{1}\neq x_{2}\in Aut\left(G\right)$ with $x_{1}E_{inn\left(G\right),Aut\left(G\right)}x_{2}$.
So define $x_{i}\left(g\right)=\left\{ \begin{array}{cc}
ag & g\notin H_{i}\\
g & g\in H_{i}\end{array}\right.$. Since $x_{i}^{2}=id$ ,$x_{i}\pi\left(g\right)x_{i}^{-1}=\pi\left(x_{i}\left(g\right)\right)=\pi\left(g\right)$
and the fact that $x_{i}\notin Inn\left(G\right)$ (because $Z\left(G\right)=Z\left(H_{i}\right)$)
it follows that \\
$tp_{qf}\left(x_{1},Inn\left(G\right),Aut\left(G\right)\right)=tp_{qf}\left(x_{2},Inn\left(G\right),Aut\left(G\right)\right)$.
Now we have to construct such a group. Notice that it is enough to
find a centerless group satisfying only the last two properties, since
we can take it's product with $\mathbb{Z}/2\mathbb{Z}$. So take $G=D_{\infty}=\left\langle a,b\left|a^{2}=b^{2}=e\right.\right\rangle $,
and $H_{a}=ker\varphi_{a}$ where $\varphi_{a}:G\to\mathbb{Z}/2\mathbb{Z}$
takes $a$ to $1$ and $b$ to $0$. In the same way we define $H_{b}$,
and finish.

\end{enumerate}
The following is the crucial claim:

\begin{claim}
\label{cla:MyLemma}Assume $G_{1}\trianglelefteq G_{2}$, $C_{G_{2}}\left(G_{1}\right)=\left\{ e\right\} $
and that $\left(G_{1},A\right)$ is a special pair. Then $\left(G_{2},A\right)$
is a special pair.
\end{claim}
\begin{proof}
First we show that $C_{G_{2}}\left(A\right)=\left\{ e\right\} $.
Suppose that $x\in C_{G_{2}}\left(A\right)$, so $xax^{-1}=a$ for
all $a\in A$. Since conjugation by $x$ (i.e. the map $h\mapsto xhx^{-1}$
in $G_{1}$) is an automorphism of $G_{1}$, (as $G_{1}$ is a normal
subgroup of $G_{2}$), it follows from $\left(G,A\right)$ being a
special pair (by remark \ref{rem:E_GA}, clause (\ref{enu:EqivDefOfE_GA}))
that it must be $id$. Hence, $x\in C_{G_{2}}\left(G_{1}\right)$,
but we assumed $C_{G_{2}}\left(G_{1}\right)=\left\{ e\right\} $ hence
$x=e$ . \\
Next assume that $xE_{G_{2},A}y$ where $x,y\in G_{2}$ and we shall
prove $x=y$. There is an isomorphism $\pi:\left\langle x,A\right\rangle _{G_{2}}\to\left\langle y,A\right\rangle _{G_{2}}$
taking $x$ to $y$ and fixing $A$. We wish to show that $x=y$,
so it is enough to show that $x^{-1}\pi\left(x\right)\in C_{G_{2}}\left(A\right)$.
This is equivalent to showing $x^{-1}\pi\left(x\right)a\pi\left(x^{-1}\right)x=a$,
i.e. $x^{-1}\pi\left(xax^{-1}\right)x=a$, i.e. $\pi\left(xax^{-1}\right)=xax^{-1}$
(remember that $\pi\left(a\right)=a$) for every $a\in A$. But $xax^{-1}$
is an element of $G_{1}$ (as $G_{1}\trianglelefteq G_{2}$), and
$\pi:\left\langle xax^{-1},A\right\rangle _{G_{1}}\to\left\langle \pi\left(xax^{-1}\right),A\right\rangle _{G_{1}}$
must be $id$ because $\left(G_{1},A\right)$ is a special pair, and
we are done.
\end{proof}
\begin{note}
If $G$ is centerless then $G\trianglelefteq Aut\left(G\right)$,
and $C_{Aut\left(G\right)}\left(G\right)=\left\{ e\right\} $.
\end{note}
\begin{conclusion}
\label{con:SpecialPairInd}Assume $G$ is centerless and $\left(G,A\right)$
is a special pair then:
\end{conclusion}
\begin{enumerate}
\item \label{enu:specialInd}$\left(G^{\alpha},A\right)$ is a special pair
for every $\alpha\in\mathbf{ord}$.
\item \label{enu:CGAin}$C_{G^{\alpha}}\left(A\right)=\left\{ e\right\} $
for every $\alpha$.
\end{enumerate}
\begin{proof}
(\ref{enu:CGAin}) follows from (\ref{enu:specialInd}). Prove (\ref{enu:specialInd})
by induction on $\alpha$. For limit ordinal, its clear from the definitions,
and for successors, the previous claim finishes the job using the
above note.
\end{proof}
\begin{conclusion}
\label{con:NormTower}Let $\gamma$ be an ordinal, $G$ a centerless
group then:
\end{conclusion}
\begin{enumerate}
\item \label{enu:CenterOfTowerIsG}$C_{G^{\gamma}}\left(G\right)=\left\{ e\right\} $.
\item $nor_{G^{\gamma}}\left(G^{\beta}\right)=G^{\beta+1}$, for $\beta<\gamma$.
\item $nor_{G^{\gamma}}^{\beta}\left(G\right)=G^{\beta}$ for $\beta\leq\gamma$.
\end{enumerate}
\begin{proof}
$\,$
\begin{enumerate}
\item Follows from conclusion \ref{con:SpecialPairInd} and from the fact
that $\left(G,G\right)$ is a special pair.
\item The direction $nor_{G^{\gamma}}\left(G^{\beta}\right)\geq G^{\beta+1}$
is clear from the definition of the action of $G^{\beta+1}$ on $G^{\beta}$.
The direction $nor_{G^{\gamma}}\left(G^{\beta}\right)\leq G^{\beta+1}$
follows from the previous clause: suppose $y\in nor_{G^{\gamma}}\left(G^{\beta}\right)$,
so conjugation by $y$ is in $Aut\left(G^{\beta}\right)$. By definition
there is $z\in G^{\beta+1}$ such that $yxy^{-1}=zxz^{-1}$ for all
$x\in G^{\beta}$, in particular - for all $x\in G$, So $y=z$ (by
(\ref{enu:CenterOfTowerIsG}))
\item By induction on $\beta$.
\end{enumerate}
\end{proof}
\begin{claim}
\label{cla:FirstBound}If $G$ is centerless and $\left(G,A\right)$
is a special pair then:
\end{claim}
\begin{enumerate}
\item \label{enu:ZFCBound}($ZFC$) $\left|G^{\alpha}\right|\leq2^{\left|A\right|}$
for all ordinals $\alpha$.
\item \label{enu:ZFBound}($ZF$) There is a one to one absolutely definable
(with parameters $G^{\alpha}$ and $A$) function from $G^{\alpha}$
into $\mathcal{P}\left({}^{\omega>}A\right)$ for each ordinal $\alpha$.
\end{enumerate}
\begin{proof}
(\ref{enu:ZFCBound}) follows from (\ref{enu:ZFBound}). The natural
way to define the function $f$ is \\
$f\left(g\right)=tp_{qf}\left(g,A,G^{\alpha}\right)$, which is a
set of equations. Luckily it is easy to encrypt equations as elements
of ${}^{\omega>}A$: We can assume that there are at least two elements
in $A$ - $a,b$ (if not, $G=\left\{ e\right\} $ because $C_{G}\left(A\right)=\left\{ e\right\} $).
Let $\sigma\left(z,\bar{a}\right)$ be a word, so it is of the form
$\ldots z^{n_{i}}a^{n_{i+1}}z^{n_{i+2}}\ldots$ where $n_{i}\in\mathbb{Z}$,
and $i=0,\ldots,m-1$. First we encrypt the exponents series with
a natural number, $m$, using the bijection $cd:{}^{\omega>}\mathbb{\omega}\to\omega$,
and then we encrypt the series of indices where $z$ appears, call
it $k$. Then we encrypt $\sigma$ by $a^{k}\mathcircumflex b\mathcircumflex a^{m}\mathcircumflex b$
and after that - the list of elements of $A$ in $\sigma$ by order
of appearance.\\
Note that our function is definable as promised.
\end{proof}
\begin{claim}
\label{cla:secondBound}If $G$ is centerless then:
\end{claim}
\begin{enumerate}
\item \label{enu:ZFCTauBound}($ZFC$) If $\left|G^{\alpha}\right|\leq\lambda$
for all ordinals $\alpha$, then $\tau_{G}<\lambda^{+}$.
\item \label{enu:ZFTauBound}($ZF$) If $\left|G^{\alpha}\right|\leq\left|A\right|$
for all ordinals $\alpha$ and a set $A$, then $\tau_{G}<\theta_{A}$.
It is enough to assume that there is a function from $A$ onto $G^{\alpha}$
for each ordinal $\alpha$.
\end{enumerate}
\begin{proof}
(\ref{enu:ZFCTauBound}) follows from (\ref{enu:ZFTauBound}), but
with choice, it is much simpler - $G_{\lambda^{+}}=\bigcup\left\{ G_{\alpha}\left|\alpha<\lambda^{+}\right.\right\} $.
Since $\left|G_{\lambda^{+}}\right|\leq\lambda$ and $\left\langle G_{\alpha}\right\rangle $
is increasing, it follows that there must be some $\alpha<\lambda^{+}$
such that $G_{\alpha}=G_{\alpha+1}$. \\
For the second part, first we show that $\tau_{G}$ is well defined.
For this we note that if $G^{\alpha}\neq G^{\alpha+1}$ then $\tau_{G^{\alpha+1},G}=\alpha+1$
(see conclusion \ref{con:NormTower}). By claim \ref{cla:ZFNormBound},
$\theta_{A}\geq\theta_{G^{\alpha+1}}>\alpha+1$. Since $\theta_{A}$
is well defined, $\tau_{G}$ is well defined as well. Applying the
same argument to $G^{\tau_{G}}$, we see that $\theta_{A}\geq\theta_{G^{\tau_{G}}}>\tau_{G}$.
\end{proof}
So as promised, we proved Thomas' theorem in a different way, without
choice:

\begin{conclusion}
\label{con:Thomas}($ZFC$) Thomas' theorem: if $G$ is a centerless
group then $\tau_{G}<\left(2^{\left|G\right|}\right)^{+}$. Moreover,
$\tau_{\kappa}<\left(2^{\kappa}\right)^{+}$.
\end{conclusion}
\begin{proof}
Taking $A=G$, $\left(G,A\right)$ is a special pair applying \ref{cla:FirstBound}
and \ref{cla:secondBound} we get the result regarding $\tau_{G}$.
Noting that $\left(2^{\kappa}\right)^{+}$ is regular and that there
are at most $2^{\kappa}$ groups of order $\kappa$ we are done.
\end{proof}
Now we deal with the case without choice. 

\begin{main theorem}
\label{cla:TauGZF}($ZF$) If $\left(G,A\right)$ is a special pair
and $G$ is a centerless group, then $\tau_{G}<\theta_{\mathcal{P}\left({}^{\omega>}A\right)}$.
\end{main theorem}
\begin{proof}
By claim \ref{cla:secondBound}, clause (\ref{enu:ZFTauBound}), we
only need to show that $\left|G^{\alpha}\right|\leq\left|\mathcal{P}\left({}^{\omega>}A\right)\right|$,
but this is exactly claim \ref{cla:FirstBound}, clause (\ref{enu:ZFBound}).
\end{proof}
Now we shall improve this by:

\begin{main theorem}
\label{thm:ZFBoundOnTauK}($ZF$) $\tau_{\left|k\right|}<\theta_{\mathcal{P}\left({}^{\omega>}k\right)}$.
\end{main theorem}
\begin{proof}
Recall that $\tau_{\left|k\right|}=\bigcup\left\{ \tau_{G}+1\left|G\right.\textrm{ is centerless and $\left|G\right|\leq\left|k\right|$}\right\} $,
but we can replace this by $\tau_{\left|k\right|}=\bigcup\left\{ \tau_{G}+1\left|G\in\mathcal{G}\right.\right\} $
where \\
$\mathcal{G}=\left\{ G\left|G\textrm{ is centerless and $G\subseteq k$}\right.\right\} $.
By the previous theorem (\ref{cla:TauGZF}) we know that $\tau_{\left|k\right|}\leq\theta_{\mathcal{P}\left({}^{\omega>}k\right)}^{\mathbf{V}'}$,
(for all $G\in\mathcal{G}$, $\left(G,G\right)$ is a special pair,
so $\tau_{G}<\theta_{\mathcal{P}\left({}^{\omega>}G\right)}\leq\theta_{\mathcal{P}\left({}^{\omega>}k\right)}$)
but we want more.\\
We may assume WLOG that $\tau_{\left|k\right|}>\omega$, since $\theta_{\mathcal{P}\left({}^{\omega>}k\right)}>\omega$
(see remark \ref{rem:Theta}). Let $\mathcal{G}'=\left\{ G\in\mathcal{G}\left|\tau_{G}\textrm{ is infinite}\right.\right\} $,
so $\tau_{\kappa}=\bigcup\left\{ \tau_{G}+1\left|G\in\mathcal{G}'\right.\right\} $.\\
 For each $G\in\mathcal{G}'$ we define a function $R_{G}:\mathcal{P}\left({}^{\omega>}k\right)\to\tau_{G}+1$
which is onto: first we define a function from $\mathcal{P}\left({}^{\omega>}k\right)$
onto $G^{\tau_{G}}$ (using claim \ref{cla:FirstBound}) , then from
$G^{\tau_{G}}$ onto $\tau_{G}$ (using claim \ref{cla:ZFNormBound},
and claim \ref{con:NormTower}), and then from $\tau_{G}$ onto $\tau_{G}+1$
(remember that $\tau_{G}\geq\omega$). \\
Let $\mathcal{B}=\left\{ \left(x,G\right)\left|G\in\mathcal{G}',x\in\mathcal{P}\left({}^{\omega>}k\right)\right.\right\} $.
Define a function $R_{1}:\mathcal{B}\to\tau_{\kappa}$ by $R\left(\left(x,G\right)\right)=R_{G}\left(x\right)$
(Note - since $R_{G}$ is definable, there is no use of $AC$). By
definition, $R_{1}$ is onto. Now it is enough to find an onto function
$R_{2}:\mathcal{P}\left({}^{\omega>}k\right)\to\mathcal{B}$. But
there is an injective function from $\mathcal{B}$ to $\mathcal{P}\left({}^{\omega>}k\right)$:
$G$ is a triple of nonempty subsets of $k$, so it is enough to know
how to encode pairs $\left(a,b\right)$ where $\emptyset\neq a\subseteq k$
and $b\subseteq{}^{\omega>}k$ as a subset $c\subseteq{}^{\omega>}k$.
For instance let $c=\left\{ x\mathcircumflex\bar{y}\left|x\in a,\bar{y}\in b\right.\right\} $.
\end{proof}
Using the following absoluteness lemma:

\begin{lem}
\label{lem:Abs}.Let $\mathbf{V}'\subseteq\mathbf{V}$ a transitive
subclass, which is a model of $ZF$. Let $\left(G,A\right)$ be a
special pair, and suppose $G,\mathcal{P}\left({}^{\omega>}A\right)\in\mathbf{V}'$.
Then, for every ordinal $\delta\in\mathbf{V}'$, the automorphism
tower $\left\langle G^{\beta}\left|\beta<\delta\right.\right\rangle $
in $\mathbf{V}'$ is the same in $\mathbf{V}$ (i.e. \\
$\mathbf{V}\models\textrm{"}\left\langle G^{\beta}\left|\beta<\delta\right.\right\rangle \textrm{ is the automorphism tower up to $\delta$}"$).
\end{lem}
Which we shall prove in the appendix, we can finally deduce:

\begin{thm}
\label{cla:MainClaimZF}$ $
\end{thm}
\begin{enumerate}
\item \label{enu:MainClaimZF1}Let $\mathbf{V}'\subseteq\mathbf{V}$ a transitive
subclass, which is a model of $ZF$. if $\mathcal{P}\left({}^{\omega>}k\right)\in\mathbf{V}'$,
then $\left(\tau_{\left|k\right|}\right)^{\mathbf{V}}=\left(\tau_{\left|k\right|}\right)^{\mathbf{V}'}<\theta_{\mathcal{P}\left({}^{\omega>}k\right)}^{\mathbf{V}'}$.
\item If $\kappa$ is a cardinal number in $\mathbf{V}'$ such that $\mathcal{P}\left(\kappa\right)\in\mathbf{V}'$,
then $\left(\tau_{\kappa}\right)^{\mathbf{V}}=\left(\tau_{\kappa}\right)^{\mathbf{V}'}<\theta_{\mathcal{P}\left(\kappa\right)}^{\mathbf{V}'}$.
\item In particular, $\tau_{\aleph_{0}}<\theta_{\mathbb{R}}^{L[\mathbb{R}]}$.
\end{enumerate}
\begin{proof}
Obviously, we need only to see (\ref{enu:MainClaimZF1}). Let \\
$\mathcal{G}=\left\{ G\left|G\textrm{ is a centerless group and }G\subseteq k\right.\right\} $.
By the assumption on $k$, it is easy to see that $\mathcal{G}^{\mathbf{V}}=\mathcal{G}^{\mathbf{V}'}$.
Hence \\
$\tau_{\left|k\right|}^{\mathbf{V}}=\bigcup\left\{ \tau_{G}+1\left|G\in\mathcal{G}^{\mathbf{V}}\right.\right\} =\bigcup\left\{ \tau_{G}+1\left|G\in\mathcal{G}^{\mathbf{V}'}\right.\right\} =\tau_{\left|k\right|}^{\mathbf{V}'}$
(the second equality is lemma \ref{lem:Abs}). By theorem \ref{thm:ZFBoundOnTauK},
we have $\tau_{\left|k\right|}^{\mathbf{V}'}<\theta_{\mathcal{P}\left({}^{\omega>}k\right)}^{\mathbf{V}'}$.
\end{proof}
If we apply lemma 1.8 from \cite{key-1}, which says that $\tau_{\kappa}^{nlg}\leq\tau_{\kappa}$
and get:

\begin{main conclusion}
\label{con:FianlConclusion}Let $\mathbf{V}'\subseteq\mathbf{V}$
be as before (but now assume $\mathbf{V}\models ZFC$) .If $\mathcal{P}\left({}^{\omega>}k\right)\in V'$,
then $\tau_{\left|k\right|}^{nlg}\leq\tau_{\left|k\right|}<\theta_{\mathcal{P}\left({}^{\omega>}k\right)}^{V'}$. 
\end{main conclusion}
\begin{note}
We actually don't need to assume that $\mathbf{V}$ is a model of
$ZFC$. $\tau_{\kappa}^{nlg}\leq\tau_{\kappa}$ is true even without
choice, and this subject will be addressed in a later work.
\end{note}

\section{\label{sec:The-descriptive-set}The descriptive set theoretic Result}

In this short section we give a descriptive set theoretic approach
into finding a bound on $\tau_{\aleph_{0}}$. We start with the definition.

\begin{defn}
\label{def:InductiveOrd}Let $M$ be structure.
\end{defn}
\begin{enumerate}
\item For a formula $\varphi\left(x,X\right)$ - a first order formula in
the language of $M$, where $x$ is a single variable and $X$ is
a mondaic variable (i.e. serve as a unary predicate - vary on subset of the structure, so not quantified
inside the formula) - we define a sequence $\left\langle X_{\alpha}^{\varphi}\subseteq M\left|\alpha\in\mathbf{Ord}\right.\right\rangle $
by:

\begin{itemize}
\item $X_{0}^{\varphi}=\emptyset$.
\item $X_{\alpha+1}^{\varphi}=X_{\alpha}\cup\left\{ x\in M\left|\varphi\left(x,X_{\alpha}^{\varphi}\right)\textrm{ is satisfied in }M\right.\right\} $.
\item $X_{\delta}^{\varphi}=\bigcup\left\{ X_{\beta}^{\varphi}\left|\beta<\delta\right.\right\} $
for $\delta$ limit.
\end{itemize}
\item For such a formula $\varphi$, let $\delta_{\varphi}=min\left\{ \alpha\left|X_{\alpha}^{\varphi}=X_{\alpha+1}^{\varphi}\right.\right\} $.
\item Let $\delta=\delta\left(M\right)$ - the inductive ordinal of the
structure - be the first ordinal such that for any such formula (allowing members of $M$ as parameters) $\varphi$,
$\delta_{\varphi}<\delta$.
\end{enumerate}
\begin{thm}
For a centerless group $G$ with set of elements $\subseteq\omega$
the height of its automorphism tower is smaller then the inductive
ordinal of the structure $\mathfrak{A}$ with universe $\omega\cup\mathcal{P}\left(\omega\right)$
the operations of $\mathbb{N}$ , membership, and $G$ (i.e. its product)
.
\end{thm}
\begin{note}
In this version of the theorem we do not need to use parameters in 
definition \ref{def:InductiveOrd}.  However the theorem holds 
even without assuming that the structure contains $G$, but
then we need parameters ($G$ can be encoded as a subset of $\omega$).
In that case this is second order number theory.
\end{note}
\begin{proof}
(sketch) By the definition it is enough to find a formula $\Delta$
such that $X_{\alpha}^{\Delta}$ encodes $G^{\alpha}$ (including
it's multiplication and inverse). By $\left(G,G\right)$ being special,
we know that we can identify members of $G^\alpha$ ($G^\alpha$ is in the automorphism tower) as
sets of finite sequences of $\omega$ (see the proof of claim \ref{cla:FirstBound}).
It is well known that the operations of $\mathbb{N}$ allow us to
encode finite sequences. Hence, much like the proof of lemma \ref{lem:Abs},
we can find a formula $\Delta'\left(x,X_{\alpha}^{\Delta}\right)$,
as in definition \ref{def:InductiveOrd}, such that $x$ satisfies
it in $\mathfrak{A}$ iff $x$ encodes a quantifier free type of an element in $G^{\alpha+1}$.
Using a similar technique we can find a formula $\Delta''\left(x,y,X_{\alpha}^{\Delta}\right)$
such that $x,y$ satisfy it iff $x,y\in G^{\alpha+1}$ and $x\circ y=id$
(i.e. the automorphism they encode). Likewise, let $\Delta'''\left(x,y,z,X_{\alpha}^{\Delta}\right)$
say that $x\circ y=z$. Now we can define $\Delta\left(x,X_{\alpha}^{\Delta}\right)$
to say that $x$ encodes a triple $\left(a,b,c\right)$ where $a\in G^{\alpha+1}$,
$b$ encodes a pair $\left(d,d^{-1}\right)$ where $d\in G^{\alpha+1}$
and $c$ encodes a triple $\left(e,f,e\circ f\right)$ where $e$
and $f$ are from $G^{\alpha+1}$. Now we have successfully encoded $G^{\alpha+1}$
as required.
\end{proof}

\section{\label{sec:BackToChoice}Back To Choice }

Applying the proof of Thomas (which used Fodor's lemma), from \cite{key-2},
we can reduce the bound on $\tau_{G}$ for some groups. The main theorem
we shall prove is:

\begin{thm}
\label{thm:Fodor}($ZFC$) Let $G$ be a centerless group and $A\subseteq G$.
If for all ordinals $\alpha$, $C_{G^{\alpha}}\left(A\right)=\left\{ e\right\} $
then $\left|G^{\alpha}\right|\leq\left(\left|G\right|^{\left|A\right|}+\aleph_{0}\right)$.
\end{thm}
Using it and claim \ref{cla:secondBound}, clause (\ref{enu:ZFCTauBound}),
we have

\begin{conclusion}
If $\left(G,A\right)$ are as in the theorem, $\tau_{G}<\left(\left|G\right|^{\left|A\right|}+\aleph_{0}\right)^{+}$
.
\end{conclusion}
We know that if $\left(G,A\right)$ is special $A$ and $G$ satisfy
the conditions of the theorem (see \ref{con:SpecialPairInd}). Hence
in particular we have:

\begin{conclusion}
If $G$ is finitely generated, then $\tau_{G}<\aleph_{1}$.
\end{conclusion}
However, we can weaken the definition of a special pair so that more
pairs $\left(G,A\right)$ will be weakly special. So we shall start
with:

\begin{defn}
\label{def:w-Special}$\,$
\begin{enumerate}
\item For a centerless group $G$, and subgroups $H_{1},H_{2}$, we say
that a homomorphism (really a monomorphism) $\varphi:H_{1}\to H_{2}$
is good if there is an automorphism $\psi:G^{\tau_{G}}\to G^{\tau_{G}}$
(so actually an inner automorphism) such that $\varphi=\psi\upharpoonright H_{1}$.
\item If $A\subseteq G$, let $E_{G,A}^{k}$ be an equivalence relation
on $G$ defined by: $xE_{G,A}^{k}y$ iff there is a good homomorphism
taking $x$ to $y$ and fixing $A$.
\item We say that the pair $\left(G,A\right)$ is weakly special if $E_{G,A}^{k}$
is $\left\{ \left(x,x\right)\left|x\in G\right.\right\} $.
\end{enumerate}
\end{defn}
\begin{rem}
If $xE_{G,A}^{k}y$ then also $xE_{G,A}y$ but not necessarily the
other direction, and so if $\left(G,A\right)$ is special, it is also
weakly special (so the name is justified)
\end{rem}
\begin{claim}
If $G_{1}$ is centerless, $G_{2}=Aut\left(G_{1}\right)$, and $\left(G_{1},A\right)$
is weakly special, then so is $\left(G_{2},A\right)$.
\end{claim}
\begin{proof}
The proof is identical to the proof of \ref{cla:MyLemma}, since conjugation
is a good homomorphism, $G_{1}\trianglelefteq G_{2}$ and $G_{1}^{\tau_{G_{1}}}=G_{2}^{\tau_{G_{2}}}$.
\end{proof}
And much like conclusion \ref{con:SpecialPairInd} we have:

\begin{conclusion}
If $\left(G,A\right)$ is (weakly) special then so is $\left(G^{\alpha},A\right)$
for every ordinal $\alpha$, and $C_{G^{\alpha}}\left(A\right)=\left\{ e\right\} $.

After giving the definition, let us prove theorem \ref{thm:Fodor}:
\end{conclusion}
\begin{proof}
Similar to the proof in \cite{key-2}.\\
Denote $\lambda=\left(\left|G\right|^{\left|A\right|}+\aleph_{0}\right)$.
The proof is by induction on $\alpha$. For $\alpha=0$ its clear.
\\
Assume $\alpha=\beta+1$. $G^{\alpha}=Aut\left(G^{\beta}\right)$
but every $\varphi\in G^{\alpha}$ is determined by $\varphi\upharpoonright A$
(because inside $G^{\alpha}$ applying $\varphi$ is the same as conjugating
by it, and because of the hypothesis). This means that $\left|G^{\alpha}\right|\leq\left|{}^{A}\left(G^{\beta}\right)\right|\leq\left(\left|G\right|^{\left|A\right|}+\aleph_{0}\right)^{\left|A\right|}=\lambda$.\\
Assume that $\alpha$ is a limit ordinal. If $\alpha<\lambda^{+}$
it is clear, so suppose that $\alpha\geq\lambda^{+}$. Assume that
$\left|G^{\alpha}\right|>\lambda$. Denote $S=\left\{ \beta<\lambda^{+}\left|cf\left(\beta\right)=cf\left(\lambda\right)\right.\right\} $.
$S$ is a stationary subset of $\lambda^{+}$. Choose $B=\left\langle h_{\beta}\left|\beta\in S\right.\right\rangle $
without repetitions such that $h_{\beta}\in G^{\beta+1}\backslash G^{\beta}$
(possible because each $G^{\beta}$ before $\alpha$ is small by the
induction hypothesis). Denote the conjugation of $A$ by $h_{\beta}$
for $\beta\in S$ by $A_{\beta}$. We know that $A_{\beta}\subseteq G^{\beta}=\bigcup\left\{ G^{\gamma}\left|\gamma<\beta\right.\right\} $
($\beta$ is a limit ordinal), so by the definition of $S$, and the
fact that $cf\left(\lambda\right)>\left|A\right|=\left|A_{\beta}\right|$,
there is $\gamma<\beta$ such that $A_{\beta}\subseteq G^{\gamma}$.
This defines a regressive function $f:S\to\lambda^{+}$ by $f\left(\beta\right)=\gamma$.
By Fodor's lemma, there is a subset $S'\subseteq S$ which is stationary
(hence $\left|S'\right|=\lambda^{+}$) such that $f\upharpoonright S'$
is constant. So there is a $\gamma$ such that for every $\beta\in S'$,
$A_{\beta}\subseteq G^{\gamma}$. But for every $\beta\in S'$, $h_{\beta}$
is determined by $h_{\beta}\upharpoonright A\in{}^{A}\left(G^{\gamma}\right)$
so $\lambda^{+}=\left|S'\right|\leq\left|{}^{A}\left(G^{\gamma}\right)\right|=\lambda$
- a contradiction.
\end{proof}

\section{Appendix}

Here we shall prove the absoluteness lemma (lemma \ref{lem:Abs}).

\begin{lem}
.Let $\mathbf{V}'\subseteq\mathbf{V}$ a transitive subclass, which
is a model of $ZF$. Let $\left(G,A\right)$ be a special pair, and
suppose $G,\mathcal{P}\left({}^{\omega>}A\right)\in\mathbf{V}'$.
Then, for every ordinal $\delta\in\mathbf{V}'$, the automorphism
tower $\left\langle G^{\beta}\left|\beta<\delta\right.\right\rangle $
in $\mathbf{V}'$ is the same in $\mathbf{V}$ (i.e. \\
$\mathbf{V}\models\textrm{"}\left\langle G^{\beta}\left|\beta<\delta\right.\right\rangle \textrm{ is the automorphism tower up to $\delta$}"$).
\end{lem}

\begin{proof}
Let $\mathfrak{T}=\left\langle G^{\beta}\left|\beta\in\mathbf{ord}^{\mathbf{V}'}\right.\right\rangle $.
We shall prove by induction on $\alpha<\delta$ that $\mathfrak{T}\upharpoonright\alpha+1$
is the automorphism tower in $\mathbf{V}$ up to $\alpha+1$. \\
For $\alpha=0$ this is clear since $G\in\mathbf{V}'$.\\
For $\alpha$ limit this follows from the definitions.\\
Suppose $\alpha=\beta+1$. By the induction hypothesis $\mathfrak{T}\upharpoonright\alpha$
is the automorphism tower in $\mathbf{V}$, so $\left(G^{\beta}\right)^{\mathbf{V}}=\left(G^{\beta}\right)^{\mathbf{V}'}\in\mathbf{V}'$.
For every $\rho\in Aut\left(G^{\beta}\right)=G^{\alpha+1}$ in $\mathbf{V}$,
we need to show that $\rho\in\left(G^{\alpha}\right)^{\mathbf{V}'}$.\\
 WLOG $A$ is a subgroup of $G$ - if not, replace it with $\left\langle A\right\rangle _{G}$
(we can define a function from $A^{<\omega}$ onto $\left\langle A\right\rangle _{G}^{<\omega}$
as in claim \ref{cla:FirstBound}). Let $\mathbf{A}=A*\left\langle x\right\rangle $
i.e. the free product of $A$ and the infinite cyclic group. As in
\ref{cla:FirstBound} there is an absolute definable function from
$A^{<\omega}$ onto $\mathbf{A}$, so $\mathcal{P}\left(\mathbf{A}\right)\in\mathbf{V}'$.
Let $\mathbf{B}=A*\left\langle x,y\right\rangle $, and by the same
reasoning $\mathcal{P}\left(\mathbf{B}\right)\in\mathbf{V}'$.\\
For every $g\in G^{\alpha}$, there is an homomorphism $\varphi_{g}$
from $\mathbf{A}$ onto $\left\langle A\cup\left\{ g\right\} \right\rangle _{G^{\alpha}}$
defined by $x\mapsto g$, and fixing $A$. By \ref{con:SpecialPairInd}
($\left(G^{\alpha},A\right)$ is special), $g\mapsto ker\left(\varphi_{g}\right)$
is injective, and absolutely definable ($ker\left(\varphi_{g}\right)$
is basically just $tp_{qf}\left(g,A,G^{\alpha}\right)$). Note that
by the induction hypothesis, $\varphi_{g}^{\mathbf{V}}=\varphi_{g}^{\mathbf{V}'}$
for $g\in G^{\beta}$. Similarly, for $g,h\in G^{\alpha}$, there
is an homomorphism $\varphi_{g,h}$ from $\mathbf{B}$ onto $\left\langle A\cup\left\{ g,h\right\} \right\rangle _{G^{\alpha}}$
fixing $A$ and taking $x$ to $g$ and $y$ to $h$, and $\left(g,h\right)\mapsto ker\left(\varphi_{g,h}\right)$
is injective.\\
The following definition allows to interpret the type of $g$ in the
type of some $\left(h_{1},h_{2}\right)$ (see example below):\\

\begin{defn}
Let $B\subseteq\mathbf{B}$ 
\end{defn}
\begin{enumerate}
\item For every $\sigma\in\mathbf{B}$, Let $\psi_{\sigma}:\mathbf{A}\to\mathbf{B}$
be the homomorphism defined by $x\mapsto\sigma$, $\psi_{\sigma}\upharpoonright A=id$.
\item For $g\in G^{\beta}$ we say that $g$ is affiliated with $B$ (denoted
$g\propto B$) if there is a word $\sigma_{g}=\sigma\left(x,y,\overline{a}\right)\in\mathbf{B}$
($\overline{a}$ are parameters from $A$) such that $ker\left(\varphi_{g}\right)=\psi_{\sigma_{g}}^{-1}\left(B\right)$.
\end{enumerate}
\begin{example}
Let $\rho\in G^{\alpha},h\in G^{\beta}$. If $B=ker\left(\varphi_{\rho,h}\right)$
then for every $g\in G^{\beta}$, $g\propto B$ iff there exists $\sigma_{g}$
such that $\varphi_{\rho,h}\left(\sigma_{g}\right)=g$ (i.e. $g\in\left\langle A\cup\left\{ h,\rho\right\} \right\rangle _{G^{\alpha}}\cap G^{\beta}$).
It could easily be verified that this is indeed true, using the equality
$\varphi_{\varphi_{\rho,h}\left(\sigma\right)}=\varphi_{\rho,h}\circ\psi_{\sigma}$
for every $\sigma\in\mathbf{B}$, and \ref{con:SpecialPairInd}. 
\end{example}
We shall find an absolute first order formula $\Delta\left(H,\mathcal{P}\left({}^{\omega>}A\right),G^{\beta}\right)$
that will say {}``$H$ is a normal subgroup of $\mathbf{A}=A*\left\langle x\right\rangle $
and there exists an automorphism $\rho\in Aut\left(G^{\beta}\right)=G^{\alpha}$
such that $H=ker\left(\varphi_{\rho}\right)$''.\\
If we succeed then if $\rho\in\left(G^{\alpha}\right)^{\mathbf{V}}$
then $\Delta\left(ker\left(\varphi_{\rho}\right),\mathcal{P}\left({}^{\omega>}A\right),G^{\beta}\right)$
will hold. Since $ker\left(\varphi_{\rho}\right)\in\mathbf{V}'$,
and $\alpha$ was absolute, there is some $\rho'\in\left(G^{\alpha}\right)^{\mathbf{V}'}$such
that $ker\left(\varphi_{\rho}\right)=ker\left(\varphi_{\rho'}\right)$
so $\rho=\rho'$ and we are done.\\
Let us describe $\Delta$. It will say that $H$ is a normal subgroup
of $\mathbf{A}$ and that for each $h\in G^{\beta}$ there exists
a subgroup $B=B_{h}\leq\mathbf{B}$ with the following properties:
\begin{enumerate}
\item $B$ is a normal subgroup of $\mathbf{B}$.
\item \label{enu:RhoIsAffiliated}$H\subseteq B$, and $B\cap\mathbf{A}=H$.
\item For every $a\in A$, $a\propto B$ and $\sigma_{a}=a$ (it follows
that $B\cap A=\left\{ e\right\} $)
\item $h\propto B$ and $\sigma_{h}=y$.
\item If $g\propto B$ and $\sigma_{1}$ and $\sigma_{2}$ witness that,
then $\sigma_{1}\sigma_{2}^{-1}\in B$.
\item If $g_{1},g_{2}\propto B$ then so is $g_{1}g_{2}$ and $\sigma_{g_{1}g_{2}}=\sigma_{g_{1}}\sigma_{g_{2}}$. 
\item \label{enu:MostImportantProp}If $g\propto B$ then there exists $g'\in G^{\beta}$
such that $g'\propto B$ and $x\sigma_{g}x^{-1}=\sigma_{g'}$.
\end{enumerate}
$B=B_{h}$ induces a monomorphism $\rho_{B}$ whose domain is $H_{B}=\left\{ g\in G^{\beta}\left|g\propto B_{h}\right.\right\} $.
It is a subgroup of $G^{\beta}$ containing $A$ and $h$ (why? because
of the conditions on $B_{h}$). Also, for every $g\in H_{B}$ define
$\rho_{B}\left(g\right)$ to be the element $g'\in G^{\beta}$ as
promised from property (\ref{enu:MostImportantProp}). In order to
show that $\rho_{B}$ is a well defined monomorphism, we note that
for every $g_{1},g_{2}\in H$, if $\sigma_{g_{1}}\sigma_{g_{2}}^{-1}\in B$
then $g_{1}=g_{2}$. \\
Why? Since $B$ is normal, $\psi_{\sigma}$ induces $\psi'_{\sigma}:\mathbf{A}\to\mathbf{B}/B$,
and so the condition $g\propto B$ becomes $ker\left(\varphi_{g}\right)=ker\left(\psi'_{\sigma_{g}}\right)$.
Now, if $\sigma_{0}\in B$, then $\psi'_{\sigma_{0}\sigma}=\psi'_{\sigma}$
so $\psi'_{\sigma_{g_{1}}}=\psi'_{\sigma_{g_{2}}}$ hence $ker\left(\varphi_{g_{1}}\right)=ker\left(\varphi_{g_{2}}\right)$
.\\
Now it an easy exercise to see that $\rho_{B}$ is as promised. \\
After defining $\rho_{B}$ we demand that for every $h_{1},h_{2}\in G^{\beta}$
and suitable $B_{1}$ and $B_{2}$, $\rho_{B_{1}}$ and $\rho_{B_{2}}$
agree on their common domain. Thus we can define $\rho_{H}=\bigcup\left\{ \rho_{B_{h}}\left|h\in G^{\beta}\right.\right\} $,
and demand that $\rho_{H}$ will be an automorphism (i.e. onto). Now
all that is left is to say that $H=ker\left(\varphi_{\rho_{H}}\right)$,
and $\Delta$ is written. \\
(There is no problem with writing this in first order, since we can
talk about finite sequences from $G^{\beta}$ using $\mathcal{P}\left({}^{\omega>}A\right)$
so we can talk about $\mathbf{B},$$\mathbf{A}$, etc).\\
Why is $\Delta$ correct? because if $\Delta\left(H,\ldots\right)$
is true, then $H=ker\left(\varphi_{\rho_{H}}\right)$ by definition.
On the other hand, if $H=ker\left(\varphi_{\rho}\right)$ for some
$\rho$, then:
\begin{itemize}
\item For each $h$, $ker\left(\varphi_{\rho,h}\right)$ will be a suitable
$B_{h}$ (by the example above). 
\item If $B_{h}$satisfy the condtions above, then $\rho_{B}\upharpoonright A=\rho\upharpoonright A$
because by condition (\ref{enu:RhoIsAffiliated}) $ker\left(\varphi_{\rho_{B}\left(a\right)}\right)=\psi_{xax^{-1}}^{-1}\left(H\right)=ker\left(\varphi_{\rho\left(a\right)}\right)$.
Hence, $\rho^{-1}\circ\rho_{B}\upharpoonright A=id$ and by \ref{con:SpecialPairInd},
$\rho_{B}\upharpoonright H_{B}=\rho\upharpoonright H_{B}$.
\end{itemize}
In conclusion, the demands on $H$ are satisfied, and we are done.
\end{proof}

\end{document}